\documentclass[11pt]{amsart}
\usepackage{amssymb}
\numberwithin{equation}{section}

\theoremstyle{definition}

\theoremstyle{remark}

\begin{document}

\vskip 1cm

\centerline {\bf Vanishing Cycles and the Inverse Problem of Potential Theory.}
\vskip.3cm
\centerline {\bf Nadya Shirokova}

\vskip 1cm

\centerline {\bf Abstract}
\vskip.3cm

   In this paper we prove the infinitesimal uniqueness theorem for the Newton potential of non simply connected bodies
   using the singularity theory approach.

 We consider the Newtonian potentials of the domains in ${\bf R}^n$
   boundaries of which are the vanishing cycles on the level hypersurface of a holomorphic function
   with isolated singularity at 0. These domains don't have to be convex, connected or 
simply connected, we also don't have any dimensional restrictions.

  We consider multiparametric families of
   such domains in the miniversal deformation of the original function. We show that
   in the parameter space any point has a neighborhood s.t. potentials of the domains corresponding
   to parameters from this neighborhood are all different as functions of the external parameter y.

\vskip 1cm

\vskip .3cm
{\bf 1.1. The inverse problem of potential theory}.

\vskip .4cm

 By the classical theorem of Newton the gravitational attraction of a solid
ball on any point outside it equals to the one of the point mass (of mass
equal to the volume of the ball) in the center of the ball [N].
 Remarkable that the converse to this statement ``if a body in ${\bf R}^3$
induces the exterior potential equal to that of a point mass, it must be a 
ball'' was proved only in 80's [Z].

 The Newton potential of body D with density $ \psi(x)$ is the
 function
  $$ I(y)=\int_D {\frac{\psi(x)dx}{{((x_1-y_1)^2+...+(x_n-y_n)^2})^{(n-2)/2}}}$$

where D is a compact domain in ${\bf R}^n$, $y\in {\bf R}^n$ is a parameter near infinity.

\vskip .2cm

 One of the questions in connection 
with this problem is whether the potential distinguishes domains.
Suppose $D_1$ and $D_2$ are homogeneous solids whose potentials agree near infinity. Must $D_1=D_2$? This is so called inverse problem of Newton potential.

  Generally the answer is negative. Such examples first appeared in geophysical literature. It was shown that there exist distinct configurations in ${\bf R}^2$,
neither of which divides the plane s.t. their potentials are the same.

 Positive answer to this question was obtained by P.S. Novikov when $D_1$ and $D_2$ are domains in ${\bf R}^3$  star-like with respect to point p internal to both domains. He proves that in this situation if $D_1$ and $D_2$ have the same external potential, they coincide [Nv].
 Local question for
the various problems in potential theory were addressed in the papers of A.I.Prilepko [P1], [P2].

\vskip .2cm

 Our approach is similar to Vassiliev's [V]. He studied Newtonian potentials of
hypersurfaces in ${\bf R}^n$ using methods of singularity theory. The ramifications
(analytic continuations) of these potentials depend on a monodromy group of the singularity. He was following the theorems on Newton and Ivory [Nw,I] which assert that a potential of a charged ellipsoid equals the constant in the interior of the ellipsoid and is constant on confocal ambient ellipsoids.

  In [V-S] the multidimensional analog of these theorems for hyperboloids in Euclidean space is found. This theorem was extended by V.Arnold to the attraction of arbitrary hyperbolic hypersurfaces.

  We address the local question, when domains are close.
Boundary of $D_{\lambda}$ is a cycle of middle dimension corresponding
to parameter $\lambda$, vanishing at the critical point on the level
surface of the deformation of an analytic function $F(x,\lambda)$.
$F(x,\lambda):({\bf C}^n\times {\bf C}^{\mu},0) \rightarrow ({\bf
C},0)$ has an isolated singular point at 0. We are restricting
ourselves to the situation when $x, \lambda$ are real. Thus our cycle
is a real hypersurface on the level set of the real analytic function.

 In this setting domains don't have to be star-like or even connected
 and may have nontrivial homologies, we also don't have any
 dimensional restrictions. Domain $D_1$ is close to $D_2$,
 i.e. corresponds to parameter ${\lambda}_1$ which is close to
 ${\lambda}_2$.

 In the next paragraphs we develop all necessary machinery from the
singularity theory to give the proof of the Main Theorem.  At first
for function f equal sum of powers and corresponding domain D we prove
infinitesimal uniqueness theorem, i.e. for a large class of
deformations called miniversal deformations we show that any parameter
p in the base has a neighborhood, s.t. all domains corresponding to parameters
in this neighborhood have different potentials.

 Next we notice that any isolated singular point can be realized in a
miniversal deformation of the sum of powers.  From the analytical
point of view this means that polynomials are dense in the
corresponding functional space, in singularity theory this follows
from the theorems of finite determinacy.

\vskip .5cm

{\bf Main Theorem}. Let $F:({\bf R}^n \times {\bf R}^{\mu},0)
\rightarrow ({\bf R},0)$ a miniversal deformation of an analytic
function with isolated singularity and a family of domains
$D_{\lambda_i}$ with boundaries - vanishing cycles on the level
surface. Then for any parameter $\lambda$ in the base of miniversal
deformation $\lambda\in \Lambda \setminus \Sigma$ ($\Sigma$ is the
discriminant) there exists a neighborhood $U_{\lambda}$ s.t. for any
${\lambda}_1,{\lambda}_2\in U_{\lambda}$, $I_{{\lambda}_1}(y)\neq
I_{{\lambda}_2}(y)$.

\vskip .2cm

 Note that since we have a multiparametric deformation, there will be
directions in the base of the miniversal deformation s.t. the volume
of the corresponding domains won't change.
 
\vskip .4cm

\vskip .7cm

{\bf 1.2. The singularity theory approach}.

\vskip .4cm

 First we introduce several definitions.

\vskip .3cm

 {\bf Definition 1.} A {\bf Milnor fibration}n is a fibration the fibers of
which are local non-singular level hypersurfaces of the functions
forming the deformation. Let $B$ be a small closed ball centered at
the origin in ${\bf C}^n$ and $c$ a generic point very close to the
origin in $C$. The corresponding manifolds
$X_{\lambda}=F^{-1}(c,\lambda) \cap B$ are called the Milnor fibers of
$F$.

\vskip .5cm

 We will consider $(n-1)$ st homology and cohomology groups of Milnor
fiber.  They are free modules over the coefficient ring and their
dimensions are equal to the multiplicity of the original critical
point. The cohomology and homology of the Milnor fibration is called
vanishing at the original critical point.

\vskip .2cm

 The fiber of Milnor fibration have a very simple topology:

\vskip .4cm

{\bf Theorem (Milnor)} [M]. The Milnor fiber $X_{\lambda}$ is homotopy
equivalent to a bouquet of spheres of the middle dimension, the number
$\mu$ of these spheres equal to the multiplicity of the initial
critical point.

\vskip .4cm

 If $f$ is a germ of a holomorphic function with a finite-multiplicity
critical point, let $\d D=\gamma_\lambda$ be a cycle vanishing at $0$
on the level surface of a function $F(x,\lambda):({\bf C}^n\times {\bf
C}^{\mu},0) \rightarrow ({\bf C},0)$ , the versal deformation of the
original function $f$, where $\lambda$ is a parameter of the
deformation.  $F(x,0)=f$.

\vskip .5cm

{\bf Definition 2}. Let $G$ be a Lie group, acting on the manifold M
and $F$ is a point of $M$. A {\bf deformation} of $f$ is a smooth map $F$
from manifold $\Lambda$ (called the base of the deformation) to $M$ at
point 0 of $\Lambda$, for which $F(0)=f$. Two deformations are
{\bf equivalent} if one can be carried to the other by the action of the
element $g(\lambda)$ of $G$ smoothly depending on the point $\lambda$
that is if $F'(\lambda)=g(\lambda)F(\lambda)$, where $g$ is a
deformation of the identity of the group. (In our case G=Diff).

\vskip .5cm

{\bf Definition 3}. Let $\varphi:(\Lambda',0) \rightarrow (\Lambda,0)$
be a smooth map. The deformation induced from $F$ by the map $\varphi$
is the deformation $\varphi^* F$ of $f$ with base $\Lambda'$, given by
the formula: $(\varphi^* F)(\lambda')=F(\varphi(\lambda))$. A
deformation $F$ is {\bf versal} if every deformation of $f$ is equivalent to
one induced from $F$. It is {\bf miniversal} if it has a minimal number of
parameters, i.e. all the parameters are essential.

\vskip .5cm

{\bf Definition 4.} The local algebra of a function $f$ at zero is the
quotient algebra of the functions by the ideal generated by the
derivatives of $f$:
$$Q_f=A_x/I_f=(\partial f/\partial x_1,...,\partial f/ \partial x_m)$$

 In [A-G-V] it was shown that one can take a versal deformation in the
form $F(x,\lambda)=f(x)+\lambda_1e_1(x)+...+\lambda_{\mu}e_{\mu}(x)$
where functions $e_k$ form a linear basis of the space $Q_f$.

\vskip .5cm

{\bf Example}.For the function equals sum of powers $f=x_1 ^N+x_2
 ^N+... x_n ^N$ one can show that the basis of the local algebra is
 formed by the monomials $e_i(x)=x_1^{i_1}....x_n^{i_n}$ where $i_k\in
 [0,N-2]$. Thus the versal deformation will have form
 $F(x,\lambda)=x_1 ^N+x_2 ^N+... x_n ^N + \lambda_1
 x_1^{N-2}...x_n^{N-2}+...+\lambda_{\mu}1$.

\vskip .2cm

 To consider the integral over the vanishing cycle we need the notion
of Gelfand-Leray form.

\vskip .5cm

{\bf Definition 5}. Let $f:{\bf R}^n \rightarrow {\bf R}$ be a smooth
function and $\omega$ be a smooth differential n-form on ${\bf
R}^n$. Smooth differential (n-1) form $\psi$ with property $df\wedge
\psi=\omega$ is called the {\bf Gelfand-Leray form} of the form $\omega$ and
denoted $\omega/df$.

\vskip .4cm

If $df$ is nonzero at some point we can show that in its neighborhood
form $\psi$ with given property exists. Restriction of this form to
any level hypersurface of a function is uniquely defined.

In our proof we will use the fact that the derivative of a function
 given as an integral of a Gelfand-Leray form over a vanishing cycle
 is expressed as an integral over the same cycle of a new form.  The
 theory of integrals over cycles vanishing in a critical point was
 developed in [AGV].

 If $\gamma_\lambda$ is a vanishing cycle on the level surface of the
holomorphic function
$F(x,\lambda)=f(x)+\lambda_1e_1(x)+...+\lambda_{\mu}e_{\mu}(x)$ then
we get the following integral representation:

 $$\partial/ {\partial}{\lambda_i}
\int_{\gamma_\lambda}\omega=\int_{\gamma_\lambda} -e_id_x\omega/d_xF$$

where $d_x\omega/d_xF$ is a Gelfand-Leray form.

\vskip .5cm

{\bf 1.3.  Potential function of the level surface}.

\vskip .5cm

 We reparametrize our problem considering instead of $x$ a new
variable $x/y$.  $D$  now is the domain close to $0$, and $y \in
S^{n-1}$ in ${\bf R}^n$

 In our setting the potential function is a function of parameter
 $\lambda \in \Lambda$ and $y \in S^{n-1}$:
 $$I_{\lambda}(y)=\int_{D_\lambda}
{\frac{\psi(x)dx}{{((x_1-y_1)^2+...+(x_n-y_n)^2})^{(n-2)/2}}}$$

where $D_\lambda$ is the domain, bounded by a vanishing cycle
$\gamma_\lambda$ and the volume form $dx=dx_1\wedge...\wedge x_n$.
Let $F$ be a smooth function in ${\bf R}^n \times {\bf R}^\mu$ and
$D_\lambda$ the hypersurface $F(x,\lambda)=0$.  Suppose that the
gradient of $F$ is nonzero at all points of the hypersurface, so that
it is smooth.

\vskip .5cm

{\bf Definition 7.} The {\bf standard charge} $\omega_F$ of the surface
$\gamma_\lambda$ is the differential form dx/dF, an (n-1)-form such
that for any tangent frame $(l_2,...l_n)$ of $\gamma_\lambda$ and a
transversal vector $l_1$ the product of the values
$\omega_\lambda(l_2,...l_n)$ and $(dF,l_1)$ is equal to the value
$dx(l_1,...l_n)$.  The natural orientation of the surface
$\gamma_\lambda$ is the orientation defined by this differential form.

\vskip .4cm

 In our case $\omega_\lambda$ is a Gelfand-Leray form restricted to
the level hypersurface and potential function has form:
 $$I_{\lambda}(y)=\int_{\gamma_\lambda}
{\frac{\psi(x)\omega_F}{{((x_1-y_1)^2+...+(x_n-y_n)^2)^{(n-2)/2}}}}$$

 We now describe the domain, for which we measure the potential. All
the most outside compact components of the level set
$\gamma_\lambda={F=0}$ get number 1. Their neighboring inside
components get number 2, etc. This way we count all the components of
$\gamma_\lambda$. Now we use Vassiliev's definition.

\vskip .5cm

{\bf Definition 8.} The {\bf Arnold cycle} of $F(x,\lambda)$ is the
manifold, oriented in such a way that the restriction to its finite
part $\gamma_\lambda$ all odd components are taken with natural
orientation while all odd components are taken with reversed
orientation.

\vskip .2cm

 Thus the potential of the domain $D_\lambda$, the boundary of which is
 $\gamma_\lambda$ is given by the integral of the Gelfand-Leray
 potential form over Arnold cycle.

\vskip .5cm

{\bf 1.4.  The Main Theorem for the sum of powers}.

\vskip .4cm

{\bf Lemma 1}. By taking linear combinations of the values of the form
$$\Omega(x,y)=((x_1-y_1)^2+...+(x_n-y_n)^2)^{-(n-2)/2}dx$$ and the
values of its derivatives over parameters $y_i$ at any given point
$y_0$ we can obtain any germ in a Taylor decomposition of the form
$\Omega(x,y)$ at 0.

\vskip .2cm

{\bf Proof}. The Taylor decomposition of $\Omega(x,y)$ at $x=0$ looks
as follows:

$$\Omega(x,y)=1+\Sigma c_{i_1,...i_n}y_1^{i_1}...y_n^{i_n}
x_1^{i_1}....x_n^{i_n}$$

 We want to be able to get any germ of the Taylor decomposition in
 finite number of steps.  Thus to have the coefficient 1 in front of
 $x_1^{i_1}...x_n^{i_n}$ it is enough to differentiate over each of
 the parameters $y_k$ $i_k$-times. After that by taking the linear
 combination of differentiated expressions with prescribed
 coefficients we can obtain a given germ.

Note that we can evaluate the resulting expression just at one point, say
$y_0=(1,0....0) \in S^n$. Thus , as we will see later, potentials can be distinguished by the moments taken at one point of the sphere.

\vskip .5cm

{\bf Theorem 1.} For the function $f=x_1^N+....+x_n^N$ all directional derivatives in
the base of versal deformation are nonzero functions.

\vskip .2cm

{\bf Proof}.  First we prove our theorem for $f = x^N + y^N + z^N$. In
this case $df = Nx^{N-1}dx + Ny^{N-1}dy + Nz^{N-1}dz$. Let $\omega^{n-2}$ be
an (n-2)-form, the differential forms which will have zero integrals
over vanishing cycle will be of the form $dfd\omega$. We want to show
that most of our forms will have nonzero integrals and we will be able
to choose a basis in the cohomology of the fiber, so that at least one
integral over the vanishing cycle, the element of the vanishing
homology, will be nonzero.

 For simplicity we will consider monomial forms. For the general (n-2)
form our algorithm will work faster, but it will be harder to trace
coefficients. 

 If we have a monomial (n-2)-form
$\omega=x_1^{a_1}x_2^{a_2}x_3^{a_3}dx_1+x_1^{b_1}x_2^{b_2}x_3^{b_3}dx_2+x_1^{c_1}x_2^{c_2}x_3^{c_3}dx_3$
forms with zero integrals will have the following presentations:

$$d\omega\wedge
df=N(x_1^{N-1}(b_3x_1^{b_1}x_2^{b_2}x_3^{b_3-1}-c_2x_1^{c_1}x_2^{c_2-1}x_3^{c_3})+$$

$$+x_2^{N-1}(a_3x_1^{a_1}x_2^{a_2}x_3^{a_3-1}-c_1x_1^{c_1-1}x_2^{c_2}x_3^{c_3})+$$

$$+x_3^{N-1}(a_2x_1^{a_1}x_2^{a_2-1}x_3^{a_3}-b_1x_1^{b_1-1}x_2^{b_2}x_3^{b_3}))dx_1dx_2dx_3$$

\vskip .5cm

 We would like to show that given a directional derivative in the base
of the versal deformation (over $\lambda$) and multiplying it by a function with any
prescribed germ (Lemma 1), i.e. taking finite linear combinations of
its derivatives over y we can obtain a basis in the vanishing cohomology of a
fiber.

 For example, if we are taking a derivative over parameter
$\lambda_\mu$ the corresponding element of the local algebra is just a
constant function 1. The integral presentation of the derivative (1.2)
will be the same as the original integral. By taking linear
combinations of the values of the derivative of $\Omega(x,\lambda)$ we can obtain any
germ of the integrated form and in particular any element of the local
cohomology ring.

  However, in general, if we differentiate over the parameter
$\lambda_i, i \ge 1$ the expression under the integral will be multiplied by
a nonconstant function $e_i$. We would like to show that the space of
forms generating cohomology ring over our cycle is the same as over
the cycle minus the set of zeros of the form which we integrate, the
directional derivative of the original function.  For the sum of
powers all elements of the local ring can be made equal
$e_1=x_1^{N-2}x_2^{N-2}x_3^{N-2}$ after multiplication by certain
monomials. Thus if we prove the theorem for the directional derivative
over $\lambda_1$, for other derivatives it will be automatic.

\vskip .2cm

 There are two operations by which we will modify our integrals in the course of the proof

 1). Assume $\lambda_\mu$ is nonzero, we observe that by
 multiplying the integral by the expressions
 $(f(x)+\lambda_1e_1(x)+...+\lambda_{\mu-1}e_{\mu-1}(x)/\lambda_\mu)^k$
 for any k, we won't change its value.
 
 2).To prove the theorem it will be enough to show that after multiplying
each basis element of the vanishing cohomology by the above expression
and then simplifying it using forms with zero integrals for some
a,b,c's we can get forms which are multiples of
$e_1=x_1^{N-2}x_2^{N-2}x_3^{N-2}$.

 Let's multiply the expression under the integral by $f^3/ \lambda_\mu$. We modify
it using forms with zero integrals. For $x_1^{3N}$ by taking $b_3=1,
b_1=2N+1$ and other a,b,c's equal zero we modify it to
$-(2N+1)x_1^{2N}x_3^{N}$. Or by taking $c_2=1, c_1=2N+1$ others zero,
to $(2N+1)x_1^{2N}x_2^{2N}$.

 Once we have monomial, which is a product of the powers of at least
two variables, we can make it divisible by
$e_1=x_1^{N-2}x_2^{N-2}x_3^{N-2}$.

 For example for $x_1^{2N}x_2^{N}$ by taking $b_3=1, b_1=N+1, b_2=N$,
a,c's are zero we get $-(N+1)x_2^Nx_2^Nx_3^N$. Or for $x_2^Nx_3^{2N}$
by taking $b_1=1,b_3=N+1, b_2=N$ one gets $(N+1)x_1^Nx_2^Nx_3^N$.  If
we will consider not just monomial, but general (n-2) forms, we could
obtain the result in 1 step.

 Next we show that for arbitrary n the proof is analogous, we have the following
``induction step'':

\vskip .5cm

{\bf Lemma 2.} Given the monomial form $\prod^k_{i=0}x_i^{\alpha_i}\cdot
 dx_1...dx_n$ we can find a form with zero integral, s.t. their
 difference will be a form $ \prod^{k+1}_{i=0}x_i^{\beta_i}\cdot
 dx_1...dx_n$.

\vskip .2cm

{\bf Proof.} Forms with zero integrals, as in the case $n=3$ will be
presented by the sum of products of $x_i^{N-1}$ and alternated sum of
monomial coefficients $a_j$ of the (n-2) form differentiated over
$x_k$, s.t. $j,k \neq i$.
 
Given
 $\omega^{n-2}=a_{1,2}(x)dx_3...dx_n+a_{1,n}(x)dx_2...dx_{n-1}+...+a_{n-1,n}dx_1...dx_{n-2}$,
 where $a_{i,j}(x)$ are the functional coefficients of the $(n-2)$ form
 containing no differentials $dx_i, dx_j$ and
 $df=Nx_1^{N-1}+...+Nx_n^{N-1}$. Then

$$df \wedge d\omega=N \sum_i x_i^{N-1}(\sum_k (-1)^k \partial
a_{i,k}/\partial x_k)$$

 After multiplying the potential form by $F(x,\lambda)^k$ we want to
kill monomials which are not the multiples of
$x_1^{N-2}.....x_n^{N-2}$. First we show the statement for
$x_1^{kN}$. Take $a_{1,2}=x_1^{kN}x_2$ and the others $a_{i,j}'s$
zero, the form with zero integral will be
$(x_1^{N-1}-x_2^Nx_1^{kN-1})dx_1...dx_n$. It's difference with $x_1^{kN}$ will have
the required form.  For any form $x_1^{\alpha_1}...x_n^{\alpha_n}
dx_1...dx_n$ take
$a_{l,l+1}=x_1^{\alpha_1-N-1}x_2^{\alpha_2}...x_l^{\alpha_l}x_{l+1}$. The
difference of the corresponding form with the original one will be
$x_1^{\alpha_1 -N-1}x_2^{\alpha_2}...x_{l+1}^Ndx_1...dx_n$.

  By adding new variables to the original form by taking sums with
forms with zero integrals like in Lemma 2, we can make any element of
the cohomology ring divisible by $e_1=x_1^{N-2}x_2^{N-2}x_3^{N-2}$.

\vskip .5cm

{\bf 1.5. The proof of the Main Theorem, Generalizations.}

\vskip .5cm

{\bf Theorem 2} . If $F(x,\lambda)$ is a versal deformation of holomorphic function
with isolated singularity, then all directional derivatives over the
parameters of the deformation are nonzero as functions of $y$.

\vskip .2cm

{\bf Proof}. In [A-G-V] it was shown that the basis of the vanishing
cohomology for the given critical point can be chosen as a part of a
basis for the sum of powers.

 We are using a fact from the singularity theory that any critical
point of multiplicity $\mu$ can be obtained in a versal deformation of
the sum of powers $x_1 ^N+x_2 ^N+... x_n ^N$ for $N>\mu +1$.

\vskip .5cm

{\bf Lemma 3}. Let the function $f:({\bf C},0) \rightarrow ({\bf
C},0)$ be holomorphic at the origin and have at the origin a critical
point of multiplicity $\mu$. Then for any $N \ge {\mu}+2$ there exists
a polynomial
$$P(x_1,....,x_n,\delta,\varepsilon_1,....,\varepsilon_n)=Q(x_1,...x_n,\delta)+
\sum^n_{j=1}(1+\varepsilon_j)x^n_j$$ possessing the properties:

\vskip .2cm

1. For fixed $\delta \neq 0,\varepsilon_1,..., \varepsilon_n$ the
function $P:({\bf C^n} \rightarrow {\bf C})$ and the function f are
equivalent in the neighborhood of the point 0.

2. $Q(x_1,...,x_n),0)=0$.

3. There exist numbers $\delta,\varepsilon_1,...,\varepsilon_n$, with
   an arbitrary small modulus, for which the hypersurface $\{x\in{\bf
   C}^n|P(x,\delta,\varepsilon)=0\}$ is nonsingular away from the
   origin.

\vskip .2cm

{\bf Proof}.  Let us take as Q a polynomial $f_{N+1}(\delta
x_1,...,\delta x_n)$, where $f_{N+1}(x_1,...,x_n)$ is the Taylor
polynomial of degree N+1 of the function f at the origin. By the
theorem of finite determinacy [A-G-V] the function in the
neighborhood of the critical point of multiplicity $\mu$ is
equivalent to its own Taylor polynomial of degree $\mu +1$.

\vskip .5cm

{\bf Corollary}. Let us consider the versal deformation of the germ of
the function$x_1^N+...+x^N_n$ at the origin. Let us denote by
$\Lambda$ the germ of the set of all values of the parameters of the
deformation for which F has a unique critical point with critical
value zero, which is equivalent to the critical point 0 of the
function f.

\vskip .4cm

 We want to show that for any critical point represented in a versal
deformation of the sum of powers we can use as the basis of the
vanishing cohomology the subset of the basis for the sum of powers.

 Since any critical point of finite multiplicity can be found obtained in the
versal deformation of the sum of powers, Milnor fibers of this critical point are included
in the Milnor fibers of the versal deformation of the sum of
powers. This embedding induces a monomorphic embedding of vanishing
homology. We can integrate the forms, generating the basis of vanishing
cohomology for the sum of powers over cycles, vanishing in a given
singular point of finite multiplicity. Theorem 2 is proved.

\vskip .5cm

{\bf Main Theorem}. Let $F:({\bf R}^n \times {\bf R}^{\mu},0)
\rightarrow ({\bf R},0)$ a miniversal deformation of an analytic
function with isolated singularity and a family of domains
$D_{\lambda_i}$ with boundaries - vanishing cycles on the level
surface. Then for any parameter $\lambda$ in the base of versal
deformation ${\lambda}\in \Lambda \setminus \Sigma$ ($\Sigma$ is the
discriminant) there exists a neighborhood $U_{\lambda}$ s.t. for any
${\lambda}_1,{\lambda}_2\in U_{\lambda}$ $I_{{\lambda}_1}(y)\neq
I_{{\lambda}_2}(y)$.

\vskip .4cm

{\bf Proof}. To prove the main theorem we just combine the results of Theorems 1 and 2. To show that
there exists a neighborhood in the parameter space in which potential
functions corresponding to the parameters from this neighborhood will be all different as functions of $y$, we use a version
of the Implicit Function Theorem. 

 We  have a functional $$I: \Lambda \rightarrow V(y)=\{f(y)\}$$ such
that for some $p \in \Lambda$, $ d_{\lambda_i} f\neq 0$, all the
directional derivatives are nonzero as functions of the parameter y.

Then on the level of the differential the map

$$d_{p} I: T_{p}\Lambda \rightarrow T_{I(p)}V$$ is injective. Thus
there exists a neighborhood in $\Lambda$ s.t.  the map $I$ is an
embedding of this neighborhood into the space of functions. This neighborhood
can be chosen as $U_{\lambda}$ from the Main Theorem and all potential functions, corresponding to parameters from $U_{\lambda}$ will be different as functions from $V(y)$ since $I$ is an embedding.

\vskip .5cm

{\bf Further directions}.

\vskip .2cm

 1. Our proof will work for any density distribution $\psi(x)$ which
has any Taylor decomposition but which is nonzero at $0$.

\vskip .2cm

2. I plan to prove the uniqueness theorem for the methaharmonic
   potential.

 \end{document}